\documentclass[12pt,leqno]{article}
\usepackage{indentfirst}

\usepackage{amsmath}
\usepackage{amscd}
\usepackage{amssymb}
\usepackage{amsthm}
\usepackage[colorlinks=true, pdfstartview=FitV, linkcolor=blue, citecolor=blue, urlcolor=blue]{hyperref}
\usepackage{color}

\def\zR{\ensuremath{\mathbb{R}}}
\def\zZ{\ensuremath{\mathbb{Z}}}

\newcommand{\Proof}{{\it Proof}}

\numberwithin{equation}{section}
\renewcommand{\theequation}{\thesection.\arabic{equation}}
\makeatletter
\newtheoremstyle{estiloobs}
  {6pt}
  {0pt}
  {\upshape}
  {}
  {}
  {}
  {.5em}
  {(\thmnumber{#2}) \thmname{\textbf{#1}}
   \thmnote{({#3})}}
\makeatother
\theoremstyle{estiloobs}
\newtheorem{remark}[equation]{Remark:}

\makeatletter
\newtheoremstyle{estilotheorem}
  {6pt}
  {0pt}
  {\slshape}
  {}
  {}
  {}
  {.5em}
  {(\thmnumber{#2}) \thmname{\textbf{#1}}
   \thmnote{({#3})}}
\makeatother
\theoremstyle{estilotheorem}

\newtheorem{lemma}[equation]{Lemma:}
\newtheorem{corollary}[equation]{Corollary:}
\newtheorem{theorem}[equation]{Theorem:}
\makeatletter
\newtheoremstyle{estilolista}
  {0pt}
  {0pt}
  {\slshape}
  {}
  {}
  {}
  {.5em}
  {(\thmnumber{#2})\thmname{\textbf{#1}}
   \thmnote{({#3})}}
\makeatother
\theoremstyle{estilolista}

\newcommand{\pedro}{\theequation}

\begin{document}

\title{Weighted inequalities for multilinear potential operators and its commutators}

\vskip 0.3 truecm

\author{Ana Bernardis, Osvaldo Gorosito and Gladis Pradolini$\,$\thanks{The authors were suportted by Consejo Nacional
de Investigaciones Cient\'\i ficas y T\'ecnicas and Universidad
Nacional del Litoral.\newline \indent Keywords and phrases:
Multilinear potential operators, commutators, weighted norm
inequalities.
}}

\date{\vspace{-0.5cm}}

\maketitle

\begin{abstract} We prove weighted strong inequalities for the
multilinear potential operator ${\cal T}_{\phi}$ and its
commutator, where the kernel $\phi$ satisfies certain growth
condition. For these operators we also obtain Fefferman-Stein type
inequalities and Coifman type estimates. On the other hand we
prove weighted weak type inequalities for the multilinear maximal
operator $\mathcal{M}_{\varphi,B}$ associated to a essentially
nondecreasing function $\varphi$ and to a submultiplicative Young
function $B$. This result allows us to obtain a weighted weak type
inequality for the operator ${\cal T}_{\phi}$.
\end{abstract}

\section{Introduction}\label{intro}

In \cite{Pe}, C. P\'erez proved weighted norm inequalities for
general potential operators. For a given nonnegative locally
integrable function $\Phi$ defined in $\zR^n$, the potential
operator $T_{\Phi}$ is defined by
\begin{equation}\label{potential}
T_{\Phi}f(x)=\int_{\zR^n}\Phi(x-y)f(y)\, dy.
\end{equation}
In a previous article E. Sawyer and R. Wheeden obtained
Fefferman-Phong type conditions on the weights and proved weighted
boundedness results for the fractional integral operator
$I_{\alpha}$ between Lebesgue spaces (see \cite{SW}). Motivated by
this paper, P\'erez considered weaker norms than those involved in
the Fefferman-Phong type conditions in \cite{SW}  and obtained
weighted boundedness results for the potential operator $T_{\Phi}$,
where the kernel $\Phi$ satisfies certain growth condition. These
norms are defined in terms of certain mapping properties of
appropriate maximal operators associated to each norm. In the same
article, the author considered the corresponding problem for the
maximal operator related to the potential operator $T_{\Phi}$
defined by
$$
M_{\tilde{\Phi}}f(x)=\sup_{x\in Q}
\frac{\tilde{\Phi}(l(Q))}{|Q|}\int_Q |f(y)|\, dy,
$$
\noindent where $l(Q)$ is the sidelength of $Q$ and
\begin{equation}\label{maxi}
\tilde{\Phi}(t)=\int_{|z|\leq t} \Phi(z)\, dz.
\end{equation}

Other authors also proved weighted norm inequalities involving the
potential operator defined in ($\ref{potential}$).

In [CUMP] the authors proved a Coifman type inequality by using
extrapolation arguments. More precisely, they showed that if $w$ is
a weight in the $A_{\infty}$ class, the $L^p(w)$ norm of $T_{\phi}$
is bounded by the $L^p(w)$ norm of $M_{\tilde{\Phi}}$. They also
obtained the corresponding result in the context of Lorentz spaces.

On the other hand, in \cite{LQY} and \cite{L1} the authors
obtained weighted $L^p$ inequalities of Fefferman-Stein type for
$T_{\Phi}$ and for the higher order commutators associated to this
operator, respectively whenever $1<p<\infty$. These estimates
involve pairs of weights of the form $(w, \mathcal{M}w)$, where
$\mathcal{M}$ are suitable maximal operators.

Two weighted norm inequalities in the spirit of those in \cite{Pe}
were proved in \cite{L2} for the higher order commutators
associated to $T_{\phi}$.

In  \cite{PW} and \cite{LYY}, weighted inequalities like the ones
described above were  proved in the general setting of the spaces
of homogeneous type.

Motivated by the work in \cite{LOPTT}, K. Moen \cite{M} considered
the multilinear fractional
 integral. For $0<\alpha <mn$ and the collection $\vec{f}=(f_1, \ldots, f_m)$ of $m$
 functions on $\mathbb{R}^n$ this operator is defined by
$$
\mathcal{I}_\alpha \vec{f}(x) =\int_{(\mathbb{R}^n)^m}
\frac{1}{(\sum_{i=1}^m|x-y_i|)^{mn-\alpha}} \prod_{i=^1}^m
f_i(y_i)\, d\vec{y}.
$$
For this operator, Moen obtained two weighted $L^p-L^q$ estimates
and a Coifman type inequality, generalizing to the multilineal
context the results given in \cite{Pe} and \cite{CUMP}. As a
consequence of the two weighted inequalities, he obtained a
generalized version of the well known characterization proved in
\cite{MW}. He also obtained weighted weak and strong inequalities
for the multi(sub)linear fractional maximal operator.

A weighted weak type inequality for $\mathcal{I}_\alpha$ and the control of the $L^p(w)$ norm of $\mathcal{I}_\alpha$ by the $L^p(Mw)$ norm of the multi(sub)linear fractional maximal operator, where $M$ is the Hardy-Littlewood maximal operator, were obtained in \cite{P}. In this article the author also obtained weighted inequalities for multi(sub)linear fractional maximal operators associated to Young functions.

On the other hand, weighted strong boundedness and weighted endpoint estimates for the multilinear
commutator of $\mathcal{I}_\alpha$ were proved in \cite{CX}. The authors also studied the multilinear fractional integral with homogeneous kernels. \\

In this paper we shall consider the multilinear potential operator defined  as

\begin{equation}
{\cal
T}_{\phi}(\vec{f})(x)=\int_{(\zR^n)^m}\phi(x-y_1,\dots,x-y_m)\prod_{i=^1}^m
f_i(y_i)\, d\vec{y},
\end{equation}
where $\phi$ is a nonnegative function defined on
$(\mathbb{R}^n)^m$. We shall also deal with the commutator
associated to this operator, given by
\begin{equation}\label{comm}
{\cal T}_{\vec{b},\phi}(\vec{f})(x)= \sum_{j=1}^m {\cal
T}_{{b_j},\phi}(\vec{f})(x),
\end{equation}
where
$$
{\cal
T}_{{b_j},\phi}(\vec{f})(x)=b_j(x){\cal T}_{\phi}(\vec{f})(x)- {\cal T}_{\phi}(f_1, \cdots , b_jf_j, \cdots , f_m)(x).
$$

The aim of this paper is to prove weighted strong and weak type inequalities
 for these multilinear operators like the ones described above. As in the case $m=1$ we
 assume that the function $\phi$ satisfies a growth condition.
More precisely, we say that a non-negative locally integrable
function $\phi$ defined in $(\zR^{n})^{m}$ satisfies a
$\mathfrak{D}$-condition or that $\phi \in \mathfrak{D}$ if there
exist two positive constants $\delta$ and $\epsilon$ such that the
inequality
\begin{equation*}\sup_{\vec{w} \in \mathcal{A}_{(2^k,1,0)}}\phi(\vec{w})\leq
\frac{C}{2^{knm}}\int_{\mathcal{A}_{(2^k,\delta,\epsilon)}}\phi(\vec{y})\, d\vec{y}
\end{equation*}
holds for every $k\in \zZ$, where  \begin{equation}\label{anillo}
\mathcal{A}_{(t,\delta,\epsilon)}=\{\vec{y}=(y_1,\cdots , y_m): \ \delta(1-\epsilon)t<\sum_{i=1}^{m}|y_i|\leq
\delta(1+\epsilon)2t\}, \quad t>0.
\end{equation}
\bigskip

Although the basic example of operators of this type is provided
by the multilinear fractional integral operator defined by the
kernel $\phi(w_1,\dots
,w_m)=\left(\sum_{i=1}^{m}|w_i|\right)^{\alpha-nm}$, for
$0<\alpha<nm$, another important example is the
multilinear Bessel potential. For $\alpha>0$ the kernel of this
operator is given by
\begin{equation*}
G_{\alpha}(x_1,\dots,x_m)=
C_{\alpha,n,m}\int_{0}^{\infty}e^{-t}e^{-\frac{\left(\sum_{i=1}^m|x_i|\right)^2}
{4t}}t^{(\alpha-nm)/2}\frac{dt}{t},
\end{equation*}
where the constant
$C_{\alpha,n,m}=\frac{1}{2^{nm}\Gamma(\alpha/2)\pi^{nm/2}}$. If
$\vec{\xi}=(\xi_1,\dots, \xi_m)$, $\xi_i\in \zR^n$, it is easy to
see that the Fourier transform of the kernel $G_{\alpha}$ valuated
in $\vec{\xi}$ is given by
\begin{equation*}
\widehat{G_{\alpha}}(\,\vec{\xi}\,)=
\left(1+4\pi^2|\,\vec{\xi}\,|\right)^{-\alpha/2}.
\end{equation*}

As in the case $m=1$ (see [G, 418]), we can prove that if $\alpha$ is a positive number and $\vec{x}=(x_1,\dots,x_m)$
with $x_i \in \zR^n$ then $G_{\alpha}(\vec{x})>0$ for every
$\vec{x}\in (\zR^n)^m$. Moreover, there exist positive constants
$C$ and $c$ that only depend on $\alpha$, $n$ and $m$ such that
\begin{equation*}
G_{\alpha}(\vec{x})\leq C e^{-|\vec{x}|/2}, \qquad {\rm when } \quad
|\vec{x}|\ge 2,
\end{equation*}
and
\begin{equation*}
c^{-1}H_{\alpha}(\vec{x})\leq G_{\alpha}(\vec{x})\leq c
H_{\alpha}(\vec{x}), \qquad {\rm for} \quad |\vec{x}|\le 2,
\end{equation*}
where $H_{\alpha}$ is a function that
satisfies
\begin{equation*}
H_{\alpha}(\vec{x})= \left\{
 \begin{array}{llllll}
 |\vec{x}|^{\alpha-nm}+1+o(|\vec{x}|^{\alpha-nm+2}) & \mbox{for $0<\alpha<nm$}, \\
 \log \frac{1}{|\vec{x}|}+1+o(|\vec{x}|^{2}) & \mbox{if $\alpha=nm$}, \\
 1+ o(|\vec{x}|^{\alpha-nm}) & \mbox{if $\alpha>nm$},
 \end{array}
 \right.
\end{equation*}
as $|\vec{x}|\rightarrow 0$. As in the linear case, $G_\alpha$ satisfies condition $\mathfrak{D}$. \\

In the setting of the spaces of homogeneous type, D. Maldonado, K.
Moen and V. Naibo \cite{MMN} defined the multilinear version of
the potential type operators studied in \cite{PW}. They proved two
weighted strong inequalities like the ones in \cite{PW} for the
case $m=1$. Observe that in the euclidean context and whenever
$K(x,y)=\phi(x-y)$, the growth condition in \cite{MMN} is more
restrictive than the condition $\mathfrak{D}$. It is easy to check
that condition $\mathfrak{D}$ is satisfied by any radial
increasing or decreasing function, and functions $\phi$
essentially constants on annuli. Since we also work with
commutator operators, our results are not included in the results
in \cite{MMN}.

\section{Preliminaries}\label{preli}

Let $X$ be a Banach function space over $\zR^n$ with respect to
the Lebesgue measure. Examples of these spaces are given by the
Lebesgue $L^p$ spaces, Lorentz spaces, Orlicz spaces and so on. Given any measurable function $f\in X$ and a cube $Q\subset
\zR^n$, we define the $X$ average of $f$ over $Q$ to be
\begin{equation*}
\|f\|_{X,Q}=\|\delta_{l(Q)}(f\chi_Q)\|_X,
\end{equation*}
where $\delta_{a}f(x)=f(ax)$ for $a>0$ and $\chi_A$ denotes the
characteristic function of the set $A$. In particular, when
$X=L^{\mathcal{B}}$, the Orlicz space associated to a Young function
$\mathcal{B}$,
\begin{equation*}
\|f\|_{X,Q}=\|f\|_{L^{\mathcal{B}},Q}=\inf\left \{\lambda>0: \frac{1}{|Q|}\int_Q
\mathcal{B}\left(\frac{|f(y)|}{\lambda}\right)\, dy\leq 1\right \}.
\end{equation*}
Particularly, when $X=L^r$, $r\ge 1$, we have that
\begin{equation*}
\|f\|_{X,Q}=||f||_{L^r,Q}=\left(\frac{1}{|Q|}\int_Q|f(y)|^r\right)^{1/r}.
\end{equation*}

Given a Banach function space $X$ there exists an associated
Banach function space $X'$ for which the generalized H\"{o}lder's
inequality,
\begin{equation*}
\int_{\zR^n}|f(x)g(x)|\, dx\leq \|f\|_X\|g\|_{{X'}},
\end{equation*}
holds. Notice that if we apply the above inequality to $\delta_{l(Q)}(f\chi_Q)$ and $\delta_{l(Q)}(g\chi_Q)$ we get that
\begin{equation*}
||fg||_{L^1,Q}\leq \|f\|_{X,Q}\|g\|_{{X'},Q}.
\end{equation*}
On the other hand, when we deal with Orlicz spaces a further
generalization of H\"older's inequality that will be useful later
is the following: If $\mathcal{A}, \mathcal{B}$ and $\mathcal{C}$
are Young functions and
$$
\mathcal{A}^{-1}(t)\mathcal{B}^{-1}(t)\le \mathcal{C}^{-1}(t)
$$
then
$$ \|f g\|_{L^\mathcal{C},Q}\le 2\|f\|_{L^\mathcal{A},
Q}\|g\|_{L^\mathcal{B}, Q},
$$
when $\mathcal{C}(t)=t$, the function $\mathcal{B}$ is called the complementary function of $\mathcal{A}$.

In general, if $\mathcal{C}$ and $\mathcal{C}_i$ are Young
functions such that
$$ \prod_{i=1}^m\mathcal{C}_i^{-1}(t) \leq \mathcal{C}^{-1}(t)
$$
then
$$
\|\prod_{i=1}^m f_i\|_{L^{\mathcal{C}},Q} \leq C \prod_{i=1}^m
\left\|f_i\right\|_{L^{\mathcal{C}_i},Q}.
$$

For more information about Orlicz spaces see \cite{O}.




\medskip

We shall say that the function $\varphi:(0,\infty)\rightarrow
(0,\infty)$ is essentially nondecreasing if there exists a
positive constant $\rho$ such that, if $t\leq s$ then
$\varphi(t)\leq \rho \varphi(s)$. We shall also suppose that
$\lim_{t\rightarrow \infty}\frac{\varphi(t)}{t}=0$. Let
$X_1,\dots,X_m$ be $m$ Banach function spaces and
$\vec{X}=(X_1,\dots,X_m)$. If $\varphi$ is a essentially
nondecreasing function, the multilinear maximal operator ${\cal
M}_{\varphi, \vec X}$ associated to $\varphi$ and $\vec X$ is
defined by
\begin{equation*}
{\cal M}_{\varphi, \vec X}\vec{f}(x)=\sup_{Q\ni
x}\varphi(|Q|)\prod_{i=1}^{m}\|f_i\|_{X_i;Q}.
\end{equation*}
If $\vec X=(X,\dots, X)$ we simply write ${\cal M}_{\varphi, X}$
and $M_{\varphi,X}$ if $m=1$. If $\varphi \equiv 1$ we write
${\cal M}_{\vec X}$, ${\cal M}_{X}$ or $M_{X}$ respectively. When
$X$ is an Orlicz space $L^{\mathcal{B}}$, sometimes we shall
denote with $||\cdot||_{\mathcal{B}}$ the norm
$||\cdot||_{L^{\mathcal{B}}}$ and with ${\cal M}_{\varphi,
\mathcal{B}}$ the operator ${\cal M}_{\varphi, L^{\mathcal{B}}}$.
In particular, when $\mathcal{B}(t)=t^p(1+\log^+t)^\alpha$, for
any $\alpha>0$ and $p\geq 1$, we shall use the notation $||\cdot||_{L^p(\log
L)^\alpha,Q}$ and ${\cal M}_{\varphi, L^p(\log L)^\alpha}$ to
denote the corresponding norm and maximal operator.
\medskip

Associated to the kernel $\phi$ we denote ${\Phi}_{\theta}$, $0<\theta\leq 1$, to be the
positive function defined for $t> 0$ by
\begin{equation}\label{tildefi}
{\Phi}_{\theta}(t)= \left[\sum_{\nu\in \zZ, \ \nu \geq \log_{1/2}
t} \left(\int_{{\cal A}(2^{-\nu},\delta,\epsilon)} \phi(\vec{y})\,
d\vec{y}\right)^{\theta}\right]^{1/ \theta}.
\end{equation}

For $\varphi(t)={\Phi}_{\theta}(t^{1/n})$ the corresponding
multilinear maximal operator ${\cal M}_{\varphi,X}$ will be denoted by
${\cal M}_{{\Phi}_{\theta},X}$.

Let $1<p_1,\dots,p_m<\infty$ such that $1/p= \sum_{i=1}^m 1/p_i$
and suppose that $M_{X_i}:L^{p_i}\rightarrow L^{p_i}$. From the
fact that ${\cal M}_{\vec{X}}\vec{f}(x)\leq
\prod_{i=1}^{m}M_{X_i}f_i(x)$ and applying H\"{o}lder's inequality
we obtain that
\begin{equation*}
{\cal M}_{\vec{X}}:L^{p_1}(\zR^n)\times\dots\times
L^{p_m}(\zR^n)\rightarrow L^p(\zR^n).
\end{equation*}

\medskip

When we deal with the commutator defined in
($\ref{comm}$) we suppose that the vector of symbols
$\vec{b}=(b_1,\dots, b_m)\in BMO^m$ and we omit this hypothesis.
Moreover, we denote $\|\vec{b}\|_*=\sup_{1\leq j\leq m} \|b_j\|_*$.

If $\mathcal{A}(t)=e^t-1$, let $||\cdot||_{\exp L,Q}$ denote the
average $||\cdot||_{L^{\mathcal{A}},Q}$. By the John Niremberg
inequality we get that $||b-b_Q||_{\exp L, Q}\leq C||b||_{*}$ for
any $BMO$ function $b$. In the proof of the results we shall use
frecuently the following inequality
$$
\frac{1}{|Q|}\int |b-b_Q||f| \leq C ||b-b_Q||_{\exp L, Q}
||f||_{L(\log L),Q}\leq C ||b||_*||f||_{L(\log L),Q}
$$
which can be proved by H\"older inequality.  \\

In this paper we shall use the following properties of the weights. We say that a weight $w$ satisfies
 a Reverse H\"older's inequality with exponent $s>1$, or equivalently, we say that $w\in RH(s)$, if there exists a positive constant $C$ such that
$$
\left(\frac{1}{|Q|}\int_Q w^s\right)^{1/s} \leq C \left(\frac{1}{|Q|} \int_Q w\right),
$$
for all cube $Q$. We say that $w\in RH_{\infty}$ if
there exists a positive constant $C$ such that the inequality
\begin{equation*}
\sup_{x\in Q} w(x)\leq \frac{C}{|Q|}\int_Q w
\end{equation*}
holds for every $Q\subset \zR^n$.\\

In all the results of this paper we shall use the following
notation: with ${\cal T}_{\vec{b}^\ell,\phi}$ we shall denote the
operator ${\cal T}_{\phi}$ when $\ell=0$ and the commutator ${\cal
T}_{\vec{b},\phi}$ when $\ell=1$. We also use $C_\ell$ to denote a
constant that depend on the norm $\|\vec{b}\|_*$ if $\ell=1$.

Sometimes we shall only estimate the term ${\cal
T}_{\vec{b_j},\phi}$ in the definition of the commutator given in
($\ref{comm}$) since the final estimate goes straightforward.

\section{Statement of the main
results}\label{section2}

\subsection{Strong type results}

In the next theorem we give Fefferman-Phong type conditions on the
weights in order to obtain the boundedness of the operator ${\cal
T}_{\phi}$ and its commutator between weighted Lebesgue spaces.

For $1\leq i,j\leq m$,
$\delta_{i,j}$ will be $1$ if $i=j$ and will be $0$ if $i\neq j$.
If $u, v_1,\dots, v_m$ are weights we write $(u,\vec{v})=(u,
v_1,\dots, v_m)$.
Let $Y_{i,j}$ and $\tilde{Y}_{i,j}$, $i,j=1,\dots,m$
be Banach function spaces. We denote ${\bf Y}$ and
$\tilde{\bf Y}$ the matrix whose components are
$Y_{i,j}$ and $\tilde{Y}_{i,j}$ respectively.

\begin{theorem}\label{acotapotencial1} Let  $\ell \in \{0,1\}$ be given. Let $1< p_1,\dots , p_m<\infty$, with $1/p=\sum_{i=1}^{m}1/p_i$, let $q$ be an
exponent satisfying $1/m<p\leq q<\infty$ and let $\phi
\in\mathfrak{D}$. For $k \in\{0,1\}$, $k\leq \ell$, and $i,j\in \{1,\cdots, m\}$, let $X^k$, $\tilde{X}^k$,
$Y_{i,j}^k$ and $\tilde{Y}_{i,j}^k$ be Banach function
spaces such that
$$
||fg||_{L(\log L)^k,Q} \leq 2 ||f||_{X^k,Q}||g||_{\tilde{X}^k,Q},
$$
$$
||fg||_{L(\log L)^{k \delta_{i,j}},Q} \leq 2
||f||_{Y_{i,j}^k,Q}||g||_{\tilde{Y}_{i,j}^k,Q},
$$
for every cube $Q$.
If $\vec{\tilde{Y}}_j^k=(\tilde{{Y}}_{1,j}^k, \dots ,
\tilde{{Y}}_{m,j}^k)$ denotes the $j$-th column in $\tilde{\bf
Y}^{k}$ suppose that for every $k\in \{0,1\}$, $k\leq \ell$,$${\cal M}_{\vec{\tilde{Y}}_j^k}:
L^{p_1}(\mathbb{R}^n)\times \dots \times L^{p_m}(\mathbb{R}^n) \to
L^{p}(\mathbb{R}^n).$$

\noindent For the weights $(u,\vec{v})$ let us write
\begin{equation}\label{B}
\mathbb{W}(\theta,\gamma,X,{\bf Y})= \max_{1\leq j\leq m}\sup_Q \
{\Phi}_{\theta}(l(Q)) \ |Q|^{1/q-1/p} \ || u^\gamma||_{X,Q}^{1/\gamma} \
\prod_{i=1}^{m}|| v_i^{-1}||_{Y_{i,j},Q}.
\end{equation}

\noindent Suppose that one of the following two conditions holds

\smallskip

 \noindent
i) $q>1$, $ M_{\tilde{X}^k}: L^{q'}(\mathbb{R}^n) \to
L^{q'}(\mathbb{R}^n)$, for every $k\in \{0,1\}$, $k\leq \ell$,
 and $(u,\vec{v})$ are weights that satisfy
\begin{equation}\label{q>11}
\max \{\mathbb{W}(1,1,X^\ell,{\bf Y}^0),\ell \ \mathbb{W}(1,1,X^0,{\bf
Y}^1)\} <\infty
\end{equation}
\noindent ii) $q\leq 1$ and $(u,\vec{v})$ are weights that satisfy
\begin{equation}\label{q<11}
\max \{\mathbb{W}(q,q,L(\log L)^{\ell q},{\bf Y}^0),\ell \
\mathbb{W}(q,1,L,{\bf Y}^1)\} <\infty
\end{equation}
Then there exists a constant $C_\ell$ such that the inequality
\begin{equation}\label{acotacion1}
\left(\int_{\zR^n} (|{\cal T}_{\vec{b}^\ell,\phi}\vec{f}|\, u)^q
\right)^{1/q}\leq
C_\ell \prod_{i=1}^{m}\left(\int_{\zR^n}(|f_i|v_i)^{p_i}\right)^{1/p_i}
\end{equation}
holds for every $\vec{f}\in L^{p_1}(v_1^{p_1})\times\dots\times
L^{p_m}(v_m^{p_m})$.
\end{theorem}
\bigskip
Notice that,  in the case $q>1$, the factor involving
the function $\phi$ in the conditions on the weights is given by
$\Phi_1(l(Q))$. Also notice that $\Phi_1(t)$ is equivalent to
$\tilde{\phi} (\delta(1+\epsilon)t)$, where $\delta$ and
$\epsilon$ are the constant in condition $\mathfrak{D}$ and
$\tilde{\phi}$ is the function defined in \eqref{maxi}. Then, when
$q>1$ we get the corresponding conditions on the weights that the
ones previously obtained in the linear cases or for the
multilinear fractional integral. On the other hand, when the
function $\phi$ is any radial increasing or decreasing function,
or $\phi$ is essentially constant on annuli, then the function
$\Phi_\theta(t)$, $0<\theta \leq 1$, is equivalent to the function
$\tilde{\phi}(\delta(1+\epsilon)t)$. In particular, when
$\phi(\vec{y})=(\sum_{i=1}^m |y_i|)^{\alpha-nm}$
we recover the results in \cite{M} for the case $\ell=0$.\\

In the following remarks we give some examples of Banach function spaces that can
be used in conditions $\eqref{q>11}$ and $\eqref{q<11}$.

\begin{remark}\label{obs1} When we consider the operator
${\cal T}_{\phi}$ ($\ell=0$), we get the theorem above for
$X^0=L^{qr}$ and $Y_{i,j}^0=L^{rp_i'}$ for every $j=1,\dots,m$ and
for some $r>1$ (with $\tilde{X}^0=L^{(qr)'}$ and
$\tilde{Y}_{i,j}^0=L^{(rp_i')'}$). Using Orlicz spaces we get a
better result. In fact, applying the examples of Young functions
give in \cite{CUP} (see page 828) we can take $X^0=L^q(\log
L)^{q-1+\delta}$ and $Y_{i,j}^0=L^{p_i'}(\log L)^{p_i'-1+\delta}$,
for every $j=1,\dots,m$ and for some $\delta>0$.
\end{remark}

\begin{remark}\label{obs2} When $\ell=1$, that is, when we consider the operator
${\cal T}_{\vec{b},\phi}$, we get the theorem above with
$X^0=L^q(\log L)^{q-1+\delta}$ , $X^1=L^q(\log L)^{2q-1+\delta}$,
$Y_{i,j}^0=L^{p_i'}(\log L)^{p_i'-1+\delta}$ for every
$j=1,\dots,m$, $Y_{i,j}^1=L^{p_i'}(\log L)^{p_i'-1+\delta}$ for
$i\neq j$ and $Y_{j,j}^1=L^{p_j'}(\log L)^{2p_j'-1+\delta}$. The
case $m=1$ of Theorem $\ref{acotapotencial1}$ improves the results
in \cite{L2}.
\end{remark}

\bigskip

As a consequence of Theorem $\ref{acotapotencial1}$ we obtain the
following Fefferman-Stein type results.

\bigskip

\begin{corollary}\label{pesos1}
Let $1< p_1,\dots, p_m<\infty$ with $1/p=\sum_{i=1}^{m}1/p_i$, $\ell \in \{0,1\}$ and let $\phi \in \mathfrak{D}$. \\

\noindent i) If $p>1$, for each $0<\delta<1$ there exists a
constant $C_{\ell,\delta}$ such that
$$
\left(\int_{\zR^n} |{\cal
T}_{\vec{b}^\ell,\phi}\vec{f}|^p\,\Pi_{i=1}^{m}u_i^{p/p_i} \right)^{1/p}\leq
C_\ell \prod_{i=1}^{m}\left(\int_{\zR^n}|f_i|^{p_i}M_{{\Phi}_1^{p},
L(\log L)^{p(1+\ell)-1+\delta}}(u_i)\right)^{1/p_i}.
$$
\noindent ii) If $p\leq 1$ and $\ell=0$, there exists a constant $C_\ell$
such that
\begin{equation*}
\left(\int_{\zR^n} |{\cal
T}_{\phi}\vec{f}|^p\,\Pi_{i=1}^{m}u_i^{p/p_i} \right)^{1/p}\leq
C_\ell\prod_{i=1}^{m}\left(\int_{\zR^n}|f_i|^{p_i}M_{{\Phi}^p_p,L}(u_i)\right)^{1/p_i}.
\end{equation*}
\noindent iii) If $p\leq 1$ and $\ell=1$ there exists a constant $C_\ell$
such that
\begin{equation*}
\left(\int_{\zR^n} |{\cal
T}_{\vec{b},\phi}\vec{f}|^p\,\Pi_{i=1}^{m}u_i^{p/p_i}
\right)^{1/p}\leq
C_\ell \prod_{i=1}^{m}\left(\int_{\zR^n}|f_i|^{p_i}M_{{\Phi}^p_p,
L^{1/p}}(u_i)\right)^{1/p_i}.
\end{equation*}
\end{corollary}

\begin{remark}
If we take $m=1$, $\ell=0$ and $\delta=[p]+1-p$ in the estimate in
$i)$ of the Corollary above we recover a result obtained in
\cite{LQY}. On the other hand, by taking $m=1$, $\ell=1$ and
$\delta=[2p]+1-2p$, the estimate in $i)$ improves the result given
in \cite{L1} for the first order commutator.
\end{remark}

\vspace{1cm}

In the following result we obtain Coifman type
inequalities for the multilinear potential operators and its
commutators.

\bigskip

\begin{theorem}\label{coifman}
Let $\phi \in \mathfrak{D}$, $\ell \in \{0,1\}$ and $w\in A_{\infty}$. Then, for every $\vec{f}=(f_1,\dots,f_m)$ with $f_i$ bounded with compact support, we get the following inequalities.\\

\noindent i) If $0<p\leq 1$ and $\ell=0$, there exists a constant $C_\ell$ such that
\begin{equation*}
\int_{\zR^n}|{\cal T}_{\phi}\vec{f}(x)|^p\, w(x)\,
dx \leq C_\ell \int_{\zR^n}[{\cal M}_{{\Phi}_p, L}\vec{f}(x)]^p\, w(x)\, dx.
\end{equation*}
ii) If $0<p\leq 1$, $\ell =1$ and $w\in RH(1/p)$, there exists a constant $C_\ell$ such that
\begin{equation*}
\int_{\zR^n}|{\cal T}_{\vec{b},\phi}\vec{f}(x)|^p\, w(x)\,
dx \leq C_\ell \int_{\zR^n}[{\cal M}_{{\Phi}_p, L(\log
L)}\vec{f}(x)]^p\, w(x)\, dx.
\end{equation*}
\noindent iii) If $p> 1$, there exists a constant $C_\ell$ such that
\begin{equation*}
\int_{\zR^n}|{\cal T}_{\vec{b}^\ell,\phi}\vec{f}(x)|^p\, w(x)\,
dx \leq C_\ell \int_{\zR^n}[{\cal M}_{{\Phi}_1, L(\log
L)^{\ell}}\vec{f}(x)]^p\, w(x)\, dx.
\end{equation*}

\end{theorem}

\bigskip

For $\ell=0$, the result above was proved in \cite{M} for
the case of the multilinear fractional integral. On the other
hand, when $m=1$ an analogous result was obtained in \cite{LYY} on
spaces of
homogeneous type.\\

The following theorem contains other weighted inequalities for the
operators ${\cal T}_{\vec{b}^\ell,\phi}$ that we shall use in the
next section to prove weak type inequalities.

\bigskip

\begin{theorem}\label{for-t-d}
Let $\phi\in \mathfrak{D}$, $\ell \in \{0,1\}$, and $u$ a weight. \\

\noindent i) If $0<p\leq 1$, there exists a constant $C_\ell$ such that
\begin{equation*}
\int_{\zR^n}|{\cal T}_{{\vec{b}}^\ell,\phi}\vec{f}(x)|^p\,
u(x)\, dx\leq C_\ell
\int_{\zR^n}[{\cal
M}_{{\Phi}_1, L(\log L)^{\ell}}\vec{f}(x)]^p\, M_{L(\log L)^{\ell}}u(x)\,
dx.
\end{equation*}
\noindent ii) If $p>1$, there exists a constant $C_\ell$ such that
\begin{equation*}
\int_{\zR^n}|{\cal T}_{\vec{b}^\ell,\phi}\vec{f}(x)|^p\, u(x)\,
dx\leq
C_\ell \int_{\zR^n}[{\cal
M}_{{\Phi}_1, L(\log L)^{\ell}}\vec{f}(x)]^p M_{L(\log L)^{[\ell
p+p]}}u(x)\, dx.
\end{equation*}
\end{theorem}

\begin{remark}
For $m=1$ and $p>1$ the result above was proved in \cite{Pe1} for
the fractional integral operator $I_{\alpha}$.
\end{remark}

\subsection{Weak type results}

The following theorem is a weighted endpoint estimate for the
multilinear maximal operator ${\cal M}_{\varphi,{\mathcal B}}$.

\begin{theorem}\label{debilmaximal} Let $\varphi$ be a essencially nondecreasing function,
$\mathcal B$ a submultiplicative Young function, $\mathcal B^m=\overbrace{\mathcal B\circ\cdots\circ \mathcal
B}^{m}$, $\psi=\mathcal
B^m\circ\varphi^{1/m}$ and $u=\prod_{i=1}^{m}u_i^{1/m}$. Then
\begin{equation*}
\sup_{\lambda>0}\frac{1}{\mathcal
B^m(1/\lambda)}u\left(\left\{x\in \zR^n:{\cal
M}_{\varphi,{\mathcal
B}}\vec{f}(x)>\lambda^m\right\}\right)^m\leq C
\prod_{i=1}^m\int_{\zR^n}\mathcal B^m(|f_i|)M_{\psi,L}u_i.
\end{equation*}
\end{theorem}

Particularly, if $0\leq\alpha<nm$ and
$\varphi(t)=t^{\alpha}$, the result above was proved in \cite{P}.
On the other hand, when $\varphi(t)=1$ and $\mathcal
B(t)=t(1+\log^+t)$ an analogous result was proved in \cite{PPTT}.
Similar estimates for several maximal operators can be found in
\cite{LOPTT} in the multilineal context, and, for example, in
\cite{FS}, \cite{Pe3} in the case $m=1$.

\bigskip

The next ``control type" result follows by applying similar
techniques to those in \cite{P} for the case of multilinear
fractional operators.

\smallskip

\begin{theorem}\label{control}
Let $\phi \in \mathfrak{D}$, $\ell \in \{0,1\}$, $\delta_0$, $\delta_1>0$ and let $u$
be a weight. Then there exists a constant $C_\ell$ such that
\begin{equation*}
\|{\cal T}_{\vec{b}^\ell,\phi}\vec{f}\|_{L^{1/m,\infty}(u)}\leq C_\ell
\|{\cal M}_{{\Phi}_1, L(\log
L)^\ell}\vec{f}\|_{L^{1/m,\infty}(M_{L(\log
L)^{\ell+\delta_{\ell}}}u)}.
\end{equation*}
\end{theorem}

\smallskip

As an easy consequence of Theorem \ref{control}, Theorem
\ref{debilmaximal} applied to the case ${\mathcal B}(t)=t$, and
the inequality $M_{L(\log L)^{\delta}} (u) \leq \prod_{i=1}^{m}
[M_{L(\log L)^{\delta}}(u_i)]^{1/m}$,
we obtain the following result for the multilinear potential operator.
\bigskip
\begin{corollary}
Let $\phi \in \mathfrak{D}$, $\delta>0$ and let
$u=\prod_{i=1}^{m}u_i^{1/m}$. Then
\begin{equation*}
\|{\cal T}_{\phi}\vec{f}\|_{L^{1/m,\infty}(u)}\leq C
\prod_{i=1}^{m}\int_{\zR^n}|f_i|M_{{\Phi}_1^{1/m},L}
(M_{L(\log L)^{\delta}}(u_i)).
\end{equation*}
\end{corollary}

\smallskip

The case $m=1$ of theorem above was proved in \cite{CPSS} for
$\phi(t)=t^{\alpha-n}$ and in \cite{P} for $m>1$ and
$\phi(t)=t^{\alpha-nm}$.

\section{Proofs of strong type results}\label{section3}


We begin this section proving the following lemma based in the discretization method developed in \cite{Pe} (see \cite{M} for the case of the multilinear fractional integral).

\begin{lemma}\label{discretizacion}
Let $f_i$, $i=1,\dots,m$, be positive functions with
compact support, let $u$ be a weight and $\phi \in \mathfrak{D}$. Let $\ell \in \{0,1\}$, $0<q\leq 1$ and, for $j\in \{1,
\dots, m\}$ let $b_j$ be a BMO function. Then there exist a constant $C_\ell$, two family of dyadic cubes
$\{Q_{k,\eta}^{0}\}_{k,\eta\in \mathbb{Z}}$ and $\{Q_{k,\eta}^{j}\}_{k,\eta\in \mathbb{Z}}$ and two family of
pairwise disjoint subsets $\{E_{k,\eta}^{0}\}$ and
$\{E_{k,\eta}^{j}\}$, $E_{k,\eta}^{0} \subset Q_{k,\eta}^{0}$
and $E_{k,\eta}^{j} \subset Q_{k,\eta}^{j}$ with
\begin{equation*}
|Q_{k,\eta}^{0}|\leq C |E_{k,\eta}^{0}| \quad \text{and} \quad
|Q_{k,\eta}^{j}|\leq C |E_{k,\eta}^{j}|
\end{equation*}
for some positive constant $C$ and for every $k,\, \eta$, such that
\begin{eqnarray*}\label{discreta}
\int_{\zR^n}[|{\cal T}_{b_j^\ell,\phi}(\vec{f})| u]^q  &\leq&
C_\ell\sum_{k,\eta}{\Phi}_q(l(Q_{k,\eta}^{0}))^q ||u^q ||_{L(\log
L)^{\ell q}, 3Q_{k,\eta}^{0}}
 \prod_{i=1}^m ||f_i||^q_{L,3Q_{k,\eta}^0} |E_{k,\eta}^0| \nonumber\\
 & & + \ell C_\ell \sum_{k,\eta}{\Phi}_q( l(Q_{k,\eta}^{j}))^q
 ||u||^q_{L,3Q_{k,\eta}^{j}}
 \prod_{i=1}^m ||f_i||^q_{L(\log L)^{\delta_{i,j}},3 Q_{k,\eta}^{j}} |E_{k,\eta}^{j}|.
\end{eqnarray*}
\end{lemma}

\bigskip

\begin{proof} Proceeding as in the proof of Proposition 3.4 in \cite{CUMP})we can assume that the weight $u$ has compact support. For each $t>0$ we set $\bar{\phi}(t)=\sup_{\vec{y}\in \mathcal{A}_{(t,1,0)}}\phi(\vec{y})$,
where $\mathcal{A}_{(t,\delta,\epsilon)}$ is the set defined in
\eqref{anillo}.
Then, if $\hat{x}-\vec{y}= (x-y_1,\dots,x-y_m)$, we get
\begin{eqnarray}\label{operador}
|{\mathcal T}_{b_j^\ell,\phi}(\vec{f})(x)|&\leq&\sum_{\nu \in
\zZ}\int_{\hat{x}-\vec{y}\in
\mathcal{A}_{(2^{\nu-1},1,0)}}|b_j(x)-b_j(y_j)|^\ell
\phi(\hat{x}-\vec{y})\prod_{i=1}^m
f_i(y_i)\, d\vec{y}\nonumber\\
&\leq&\sum_{\nu \in
\zZ}\bar{\phi}(2^{\nu-1})\int_{\sum_{i=1}^{m}|x-y_i|<2^{\nu}}|b_j(x)-b_j(y_j)|^\ell
\prod_{i=1}^m
f_i(y_i)\, d\vec{y}\nonumber\\
&\leq&\sum_{\nu \in \zZ}\sum_{\{Q:l(Q)=2^{\nu}\}}\bar{\phi}\left(\frac{
l(Q)}{2}\right)\chi_Q(x) \int_{(3Q)^m}|b_j(x)-b_j(y_j)|^\ell
\prod_{i=1}^m f_i(y_i)\, d\vec{y} \nonumber\\\\
&\leq&\sum_{Q\in {\cal D}}\bar{\phi}\left(\frac{
l(Q)}{2}\right)|b_j(x)-(b_j)_Q|^\ell\chi_Q(x) \prod_{i=1}^m\int_{3Q} f_i(y_i)\, dy_i
\nonumber\\
&+& \ell
\sum_{Q\in {\cal D}}\bar{\phi}\left(\frac{
l(Q)}{2}\right)\chi_Q(x) \prod_{i=1, \ i\neq j}^m\left(\int_{3Q} f_i(y_i)\, dy_i\right) \nonumber\\
& & \times \left(\int_{3Q} |b_j(y_j)-(b_j)_Q|^\ell f_j(y_j)\, dy_j\right).
\nonumber
\end{eqnarray}
From the above estimate, using the generalized H\"older inequality and the estimate $||(b-b_Q)^q||_{exp L^{1/q},Q} \leq C||b||_{*}^q$ that holds for any BMO function $b$, we get
\begin{eqnarray*}
\int_{\zR^n}\left[|{\cal T}_{b_j^\ell,\phi}(\vec{f})|u\right]^q
&\leq& \sum_{Q\in {\cal D}}\bar{\phi}\left(\frac{
l(Q)}{2}\right)^q\left(\int_Q |b_j-(b_j)_Q|^{\ell
q}u^q\right)|3Q|^{mq} \prod_{i=1}^m ||f_i||_{L,3Q}^q
\nonumber\\
&+&  \ell
\sum_{Q\in {\cal D}}\bar{\phi}\left(\frac{
l(Q)}{2}\right)^q\left(\int_Q u^q\right)\prod_{i\neq j} \left(\int_{3Q}
f_i\right)^q \left(\int_{3Q}
|b_j - (b_j)_Q|^\ell f_j\right)^q
\nonumber\\
&\leq& C ||b_j||_*^{\ell q} \sum_{Q\in {\cal D}}\bar{\phi}\left(\frac{
l(Q)}{2}\right)^q ||u^q||_{L(\log L)^{\ell q},Q}|3Q|^{mq+1} \prod_{i=1}^m ||f_i||_{L,3Q}^q
\nonumber\\
&+& \ell C ||b_j||_*^{\ell q}
\sum_{Q\in {\cal D}}\bar{\phi}\left(\frac{
l(Q)}{2}\right)^q ||u^q||_{L,Q}|3Q|^{mq+1}\prod_{i=1}^m
||f_i||_{L(\log L)^{\delta_{i,j}},3Q}^q\\
&=& C ||b_j||_*^{\ell q}\sum_{Q\in {\cal D}}\bar{\phi}\left(\frac{
l(Q)}{2}\right)^q ||u^q||_{L(\log L)^{\ell q},Q}|3Q|^{mq+1} \prod_{i=1}^m ||f_i||_{L,3Q}^q
\nonumber\\
&+& \ell C ||b_j||_*^{\ell q}
\sum_{Q\in {\cal D}}\bar{\phi}\left(\frac{
l(Q)}{2}\right)^q ||u||_{L,Q}^q |3Q|^{mq+1}\prod_{i=1}^m ||f_i||_{L(\log L)^{\delta_{i,j}},3Q}^q.
\end{eqnarray*}
Notice that, since $q\leq 1$ in the last line of the inequalities
above we have used that $||u^q||_{L,Q} \leq ||u||_{L,Q}^q$.

From the results in \cite{M} we get the following Calder\'on-Zygmund decomposition. For $\vec{h}=(h_1,\dots,h_m)$ let
$$
{\cal M}_{3{\cal D}}\vec{h}(x)=\sup_{x\in Q\in {\cal
D}}\prod_{i=1}^m||h_i||_{L,3Q}
$$
and
$$
\mathcal{D}_k=\{x\in \mathbb{R}^n: \ {\cal M}_{3{\cal
D}}\vec{h}(x)>a^k\}
$$
with $a>6^n ||\mathcal{M}_L||$, where $||\mathcal{M}_L||$ is the
constant from the $L^1 \times \dots \times L^1 \to L^{1/m,\infty}$
inequality for $\mathcal{M}_L$. It was proved in \cite{M} that
there exists a family $Q_{k,\eta}$ of maximal, disjoint and dyadic
cubes such that ${\cal D}_k=\bigcup_{\eta \in
\mathbb{Z}}Q_{k,\eta},$ and
\begin{equation*}
a^k<\prod_{i=1}^m||h_i||_{L,3Q_{k,\eta}}\leq 2^{nm}a^k.
\end{equation*}
Moreover, if $E_{k,\eta}=Q_{k,\eta}\setminus {\cal D}_{k+1}$ then
$\{E_{k,\eta}\}_{k,\eta}$ is a disjoint family of sets that satisfies
$|Q_{k,\eta}|\leq C |E_{k,\eta}|$, for some positive constant $C$.

Observe that the cubes $Q_{k,\eta}$ and the sets $E_{k,\eta}$ depend on the function $\vec{h}$. We shall apply the above
results to two different functions. In one case we shall use $\vec{h}=\vec{f}$ and in this case
we shall denote with $Q_{k,\eta}^{0}$ and $E_{k,\eta}^{0}$ the
corresponding cubes and sets. In the other case we shall apply the above decomposition with
 $\vec{h}=(h_1^{j}, \dots, h_m^{j})$, where $h_i^{j}=f_i$ if $i\neq j$ and $h_j^{j}=u$.
  In this second case we shall denote with $Q_{k,\eta}^{j}$ and $E_{k,\eta}^{j}$ the corresponding cubes and sets. On the other hand, let
\begin{equation*}
{\cal F}_k^{j}=\{Q\in {\cal D}:
a^k<\prod_{i=1}^m ||h_i^{j}||_{L,3Q} \leq
a^{k+1}\}
\end{equation*}
and let ${\cal F}_k^{0}$ be the above set with $f_i$ instead of $h_i^{j}$.
Then
\begin{eqnarray*}
\int_{\zR^n}\left[|{\mathcal T}_{b_j^\ell,\phi}(\vec{f})|
u\right]^q &\leq& C_\ell \sum_{k\in \zZ} a^{(k+1)q}\sum_{Q\in
{\cal F}_k^{0}}\bar{\phi}\left(\frac{l(Q)}{2}\right)^q|3Q|^{mq+1}
||u^q||_{L(\log L)^{\ell q}, Q}\nonumber\\
& & + \ell C_\ell \sum_{k\in
\zZ} a^{(k+1)q}\sum_{Q\in {\cal F}_k^{j}}\bar{\phi}\left(\frac{
l(Q)}{2}\right)^q|3Q|^{mq+1}
||f_j||_{L(\log L),3Q}^q\nonumber\\
&\leq& C_{\ell}
\sum_{k\in
\zZ}a^{(k+1)q}\sum_{\eta}\sum_{Q\in {\cal F}_k^{0}: Q\subset Q_{k,\eta}^{0}}\bar{\phi}\left(\frac{
l(Q)}{2}\right)^q|3Q|^{mq+1}
||u^q||_{L(\log L)^{\ell q},3Q}\nonumber\\
& & + \ell C_\ell \sum_{k\in
\zZ}a^{(k+1)q}\sum_{\eta}\sum_{Q\in {\cal F}_k^{j}: Q\subset Q_{k,\eta}^{j}}\bar{\phi}\left(\frac{
l(Q)}{2}\right)^q|3Q|^{mq+1}
||f_j^q||_{L^{1/q}(\log L),3Q}\nonumber\\
\end{eqnarray*}

It is easy to prove (see for example \cite{L2}) that for any dyadic cube $Q_0$ with $l(Q_0)=2^{-\nu_0}$ and any Young function $\psi$, there exists a
positive constant $C$ such that
$$\sum_ {\{Q\subset Q_0, l(Q)=2^{-\nu}\}} |3 Q| \ ||f||_{L^\psi, 3 Q} \leq C |3 Q_0| \ ||f||_{L^\psi,3 Q_0}.
$$
On the other hand, since $\phi \in \mathfrak{D}$,
\begin{eqnarray*}
\sum_{\nu\geq \nu_0} [2^{-\nu n m} \bar{\phi}(2^{-\nu-1})]^q
&\leq& C \sum_{\nu\geq
\nu_0}\left(\int_{\mathcal{A}_{(2^{-\nu-1},\delta, \epsilon)}}
\phi(\vec{y})\, d\vec{y}\right)^q\\
&=& C {\Phi}_q(l(Q_0))^q.\\
\end{eqnarray*}

\noindent From the results above we get that there exists $C>0$
such that

\begin{eqnarray*}
\sum_{\{Q: Q\subset Q_0\}}\bar{\phi}\left(\frac{l(Q)}{2}\right)^q|3 Q|^{mq+1} ||f||_{L^\psi, 3 Q}
&=& 3^{nmq} \sum_{\nu\geq
\nu_0} [2^{-\nu n m} \bar{\phi}(2^{-\nu-1})]^q \\
& & \times \sum_ {\{Q\subset Q_0, l(Q)=2^{-\nu}\}} |3 Q|
||f||_{L^\psi, 3 Q}\\
&\leq& C {\Phi}_q( l(Q_0))^q|3 Q_0|||f||_{L^\psi,3 Q_0}
\end{eqnarray*}
where $C$ depends on $\delta$, $\epsilon$ in condition
$\mathfrak{D}$. Using
 the inequality above with $\psi(t)=t(1+\log^+t)^{\ell q}$ and  $\psi(t)=t^{1/q}(1+\log^+t)$ we get
\begin{eqnarray*}
\int_{\zR^n}\left[|{\mathcal T}_{b_j^\ell,\phi}(\vec{f})|
u\right]^q &\leq& C_{\ell}\sum_{k,\eta}{\Phi}_q(
l(Q_{k,\eta}^{0}))^q||u^q||_{L(\log L)^{\ell q},3Q_{k,\eta}^{0}}
 \prod_{i=1}^m ||f_i||^q_{L,3 Q_{k,\eta}^{0}} |E_{k,\eta}^{0}|\\
 & & + \ell C_\ell \sum_{k,\eta}
 {\Phi}_q( l(Q_{k,\eta}^{j}))^q||u||^q_{L, 3 Q_{k,\eta}^{j}}
 \prod_{i=1}^m ||f_i||^q_{L(\log L)^{\delta_{i,j}}, 3Q_{k,\eta}^{j}} |E_{k,\eta}^{j}|.
\end{eqnarray*}

\end{proof}

\bigskip

\noindent {\it Proof of theorem $\ref{acotapotencial1}$}: Let us
first consider the case $q>1$. It is enough to prove that, for
every $j=1, \dots , m$,
\begin{equation*}
\int_{\zR^n}|{\cal T}_{b_j^\ell,\phi}(\vec{f})(x)|u(x)g(x)\,
dx\leq C_{\ell}\|g\|_{q'}\prod_{i=1}^{m}\|f_iv_i\|_{p_i},
\end{equation*}
for all $g\in
L^{q'}(\zR^n)$, with $g\geq 0$, and for all positive function $f_i$ bounded with compact support.
From Lemma $\ref{discretizacion}$ with $q=1$, the generalized H\"{o}lder's
inequality and condition ($\ref{q>11}$)
we have that
\begin{eqnarray*}
\int_{\zR^n}|{\cal T}_{b_j^\ell,\phi}\vec{f}|u g
&\leq&C_{\ell}\sum_{k,\eta}{\Phi}_1(3
l(Q_{k,\eta}^{0}))||ug||_{L(\log L)^\ell,
Q_{k,\eta}^{0}}\prod_{i=1}^m ||f_i||_{L,3 Q_{k,\eta}^{0}}
|E_{k,\eta}^{0}|\nonumber\\
& & + \ell C_\ell \sum_{k,\eta}{\Phi}_1(3
l(Q_{k,\eta}^{j}))||ug||_{L, 3 Q_{k,\eta}^{j}}\prod_{i=1}^m ||f_i||_{L(\log L)^{\delta_{i,j}},3 Q_{k,\eta}^{j}}
|E_{k,\eta}^{j}|\nonumber\\
&\leq&C_{\ell}\sum_{k,\eta} ||g||_{\tilde{X}^\ell, 3 Q_{k,\eta}^{0} }\prod_{i=1}^m ||f_i v_i||_{\tilde{Y}_{i,j}^0,3 Q_{k,\eta}^{0}}
|E_{k,\eta}^{0}|^{1/p+1/q'}\nonumber\\
& & + \ell C_\ell \sum_{k,\eta} ||g||_{\tilde{X}^0, 3 Q_{k,\eta}^{j}}\prod_{i=1}^m ||f_i v_i||_{\tilde{Y}_{i,j}^1,3 Q_{k,\eta}^{j}}
|E_{k,\eta}^{j}|^{1/p+1/q'}\nonumber\\
\end{eqnarray*}

By H\"{o}lder's inequality and the fact that $p\le q$ we obtain
that
\begin{eqnarray*}
\int_{\zR^n}|{\cal T}_{b_j^\ell\phi}\vec{f}| u g&\leq&
C_{\ell}\left(\sum_{k,\eta} ||g||_{\tilde{X}^\ell, 3
Q_{k,\eta}^{0}}^{q'}
|E_{k,\eta}^{0}|\right)^{1/q'}\left(\sum_{k,\ell}
|E_{k,\eta}^{0}|^{q/p} \prod_{i=1}^m ||f_i
v_i||_{\tilde{Y}_{i,j}^0,3 Q_{k,\eta}^{0}}^q \right)^{1/q}
\nonumber\\
& & + \ell C_\ell \left(\sum_{k,\ell} ||g||_{\tilde{X}^0, 3 Q_{k,\eta}^{j}}^{q'}
|E_{k,\eta}^{j}|\right)^{1/q'}\left(\sum_{k,\ell}|E_{k,\eta}^{j}|^{q/p}\prod_{i=1}^m ||f_i v_i||_{\tilde{Y}_{i,j}^1,3 Q_{k,\eta}^{j}}^q \right)^{1/q}
\nonumber\\
&\leq&
C_{\ell}\left(\int_{\mathbb{R}^n}[M_{\tilde{X}^\ell}(g)]^{q'}\right)^{1/q'}
\left(\sum_{k,\eta} |E_{k,\eta}^{0}|
\prod_{i=1}^m ||f_i v_i||_{\tilde{Y}_{i,j}^0,3 Q_{k,\eta}^{0}}^p \right)^{1/p}
\nonumber\\
& & + \ell
C_\ell \left(\int_{\mathbb{R}^n}[M_{\tilde{X}^0}(g)]^{q'}\right)^{1/q'}
\left(\sum_{k,\eta} |E_{k,\eta}^{j}|
\prod_{i=1}^m ||f_i v_i||_{\tilde{Y}_{i,j}^1,3 Q_{k,\eta}^{j}}^p \right)^{1/p}
\nonumber\\
&\leq&
C_{\ell}\left(\int_{\zR^n}M_{\tilde{X}^\ell}(g)^{q'}\right)^{1/q'}
\left(\int_{\zR^n}
\left[\mathcal{M}_{\vec{\tilde{Y}}_{j}^0}(\vec{fv})\right]^{p}
\right)^{1/p}
\nonumber\\
& & + \ell C_\ell
\left(\int_{\zR^n}M_{\tilde{X}^0}(g)^{q'}\right)^{1/q'}\left(\int_{\zR^n}
\left[\mathcal{{M}}_{\vec{\tilde{Y}}_{j}^1}(\vec{fv})\right]^{p}
\right)^{1/p},
\end{eqnarray*}
where $\vec{fv}=(f_1 v_1, \dots , f_m v_m)$. Then by the
hypotheses on the boundedness properties of
$\mathcal{M}_{\vec{\tilde{Y}}_{j}^k}$ and $M_{\tilde{X}^k}$,
$k\in \{0,1\}$, $k\leq \ell$, we conclude the proof of the case $q>1$.

Suppose now that $1/m<p\leq q\leq 1$. From Lemma \ref{discretizacion}, by the generalized H\"{o}lder's
inequality, condition ($\ref{q<11}$) and the fact that $p\leq q$  we
obtain that
\begin{eqnarray*}
\left(\int_{\zR^n}\left[|{\mathcal T}_{b_j^\ell,\phi}(\vec{f})(x)|
\, u(x)\right]^q\, dx\right)^{1/q} &\leq&
C_{\ell}\left(\sum_{k,\eta} \prod_{i=1}^m ||f_i
v_i||^q_{\tilde{Y}_{i,j}^0,3 Q_{k,\eta}^{j,0}}
|E_{k,\eta}^{j,0}|^{q/p}\right)^{1/q}\nonumber\\
& & + \ell C_\ell \left(\sum_{k,\eta} \prod_{i=1}^m ||f_i v_i||^q_{\tilde{Y}_{i,j}^1,3 Q_{k,\eta}^{j}}
|E_{k,\eta}^{j}|^{q/p}\right)^{1/q}\nonumber\\
&\leq& C_\ell \left(\int_{\mathbb{R}^n}
[\mathcal{M}_{\vec{\tilde{Y}}_{j}^0}(\overrightarrow{fv})]^p
\right)^{1/p}
\nonumber\\
& & +  \ell C_\ell \left(\int_{\mathbb{R}^n}
[\mathcal{M}_{\vec{\tilde{Y}}_{j}^1}(\overrightarrow{fv})]^p
\right)^{1/p},\\
\end{eqnarray*}
and from the hypotheses on the boundedness of the multilinear
maximal operator $\mathcal{M}_{\vec{\tilde{Y}}_{j}^k}$, $k\in \{0,1\}$, $k\leq \ell$, we
finally get the theorem. \ $\square$

\bigskip

\noindent {\it Proof of Corollary \ref{pesos1}}: For the case
$p>1$ and $\ell=0$, the proof follows from Theorem
\ref{acotapotencial1} and Remark \ref{obs1} by taking $p=q$,  $X^0=L^p(\log
L)^{p-1+\delta}$ for some $0<\delta<1$, $Y_{i,j}^0$ any Orlicz
space such that the maximal operator associated to the
space  $\tilde{Y}_{i,j}^0$ is bounded from $L^{p_i}$ into $L^{p_i}$,
$i=1,\dots, m$, $u=\Pi_{i=1}^{m}u_i^{1/p_i}$ and
$v_i=[M_{{\Phi}_1^{p}, L(\log L)^{p-1+\delta}}(u_i)]^{1/p_i}$.
In fact, we use the generalized H\"older's inequality with
$\mathcal{C}(t)=t^p(1+\log^+ t)^{p-1+\delta}$ and
$\mathcal{C}_i(t)=t^{p_i}(1+\log^+ t)^{p-1+\delta}$, $1\leq i\leq m$
to estimate $||u||_{X^0}$ and the inequality
$M_{{\Phi}_1^{p},L(\log L)^{p-1+\delta}}(u_i)(x) \geq
{\Phi}_1(l(Q))^{p}||u_i||_{L(\log L)^{p-1+\delta},Q}$, for all
$x\in Q$, to estimate $||v_i^{-1}||_{Y_{i,j}^0}$. In a similar way
we can prove the result for $p\leq 1$ and $\ell=0$.\\
For the case $p>1$ and $\ell=1$, as in the previous case, the
proof follows from Theorem \ref{acotapotencial1} and Remark
\ref{obs2} by taking $p=q$, $Y_{i,j}^0$ and ${Y}_{i,j}^1$ be any
Orlicz spaces such that the maximal operators associated to the
spaces $\tilde{Y}_{i,j}^0$ and $\tilde{Y}_{i,j}^1$ are bounded
from $L^{p_i}$ into $L^{p_i}$, $i=1,\dots, m$, $X^0=L^p(\log
L)^{p-1+\delta}$, $X^1=L^p(\log L)^{2p-1+\delta}$ and the weights
$u=\Pi_{i=1}^{m}u_i^{1/p_i}$ and $v_i=[M_{{\Phi}_1^{p}, L(\log
L)^{2p-1+\delta}}(u_i)]^{1/p_i}$. For the case $p\leq1$ and
$\ell=1$ take $u=\Pi_{i=1}^{m}u_i^{1/p_i}$ and
$v_i=[M_{{\Phi}^p_p, L^{1/p}(\log L)^{p}}(u_i)]^{1/p_i}$.
\bigskip

\noindent {\it Proof of theorem $\ref{coifman}$} :
We start proving the case $p\leq 1$.
By taking $q=p$ and $u=w^{1/p}$ in Lemma \ref{discretizacion} we obtain
\begin{eqnarray*}
\int_{\zR^n}|{\cal T}_{b_j^\ell,\phi}(\vec{f})|^p \ w  &\leq&
C_\ell\sum_{k,\eta}{\Phi}_p( l(Q_{k,\eta}^{0}))^p ||w||_{L(\log
L)^{\ell p}, 3Q_{k,\eta}^{0}}
\prod_{i=1}^m ||f_i||^p_{L,3 Q_{k,\eta}^0} |E_{k,\eta}^0|\\
& & + \ell C_\ell \sum_{k,\eta}{\Phi}_p(
l(Q_{k,\eta}^{j}))^p||w^{1/p}||^p_{L,3 Q_{k,\eta}^{j}}
\prod_{i=1}^m ||f_i||^p_{L(\log L)^{\delta_{i,j}},3
Q_{k,\eta}^{j}} |E_{k,\eta}^{j}|.
\end{eqnarray*}

Since $w\in A_{\infty}$ then $w$ satisfies the reverse
H\"{o}lder's condition $RH(s)$ for some $s>1$. Then, from the fact that $\|.\|_{L(\log L)^{\ell p}, Q}\leq \|.\|_{L^s, Q}$ for any $s>1$, the
condition $RH(s)$ and the condition $RH(1/p)$ for the case $\ell=1$, we get that
\begin{eqnarray*}
\int_{\zR^n}|{\cal T}_{b_j^\ell,\phi}(\vec{f})|^p \ w  &\leq&
C_\ell\sum_{k,\eta}{\Phi}_p( l(Q_{k,\eta}^{0}))^p ||w||_{L,
3Q_{k,\eta}^{0}}
 \prod_{i=1}^m ||f_i||^p_{L,3 Q_{k,\eta}^0} |E_{k,\eta}^0|\\
 & & +  \ell C_\ell \sum_{k,\eta}{\Phi}_p( l(Q_{k,\eta}^{j}))^p
 ||w||_{L,3 Q_{k,\eta}^{j}}
 \prod_{i=1}^m ||f_i||^p_{L(\log L)^{\delta_{i,j}},3 Q_{k,\eta}^{j}} |E_{k,\eta}^{j}|.
\end{eqnarray*}
Since $w\in A_{\infty}$  and $|Q_{k,\eta}|\leq C|E_{k,\eta}|$ we obtain
that $w(Q_{k,\eta})\leq C w(E_{k,\eta})$, so that
\begin{eqnarray*}
\int_{\zR^n}|{\cal T}_{b_j^\ell,\phi}(\vec{f})|^p \ w  &\leq&
C_\ell \sum_{k,\eta}{\Phi}_p( l(Q_{k,\eta}^{0}))^p
 \prod_{i=1}^m ||f_i||^p_{L,3 Q_{k,\eta}^0} w(E_{k,\eta}^0)\\
 & & + \ell C_\ell \sum_{k,\eta}{\Phi}_p(l(Q_{k,\eta}^{j}))^p
  \prod_{i=1}^m ||f_i||^p_{L\log L,3 Q_{k,\eta}^{j}}
  w(E_{k,\eta}^{j})\\
&\leq& C_\ell \sum_{k,\eta}\int_{E_{k,\eta}^0}[{\cal
M}_{{\Phi}_p}\vec{f}(x)]^p w(x) \, dx\\
 & & + \ell C_\ell \sum_{k,\eta}\int_{E_{k,\eta}^{j}}[{\cal
M}_{{\Phi}_p,L\log L}\vec{f}(x) w(x)]^p \, dx\\
&\leq & C_\ell \int_{\zR^n}[{\cal
M}_{{\Phi}_p,L(\log L)^{\ell}}\vec{f}(x)]^p w(x) \, dx.\\
\end{eqnarray*}
The case $p>1$ follows from the case $p=1$ by applying
the extrapolation theorem in \cite{CUMP}. $\square$

\bigskip

We shall use the next two results to prove Theorem $\ref{for-t-d}$. The first one is a
corollary of Lemma \ref{discretizacion} and the second one was proved in \cite{CN}.

\begin{corollary}\label{TM} Let $\phi\in \mathfrak{D}$ and let $v$ be a weight
satisfying the $RH_{\infty}$ condition. If $u$ is a weight, then
there exists a positive constant $C_\ell$ such that
\begin{eqnarray*}
\int_{\zR^n}|{\cal T}_{\vec{b}^\ell,\phi}\vec{f}(x)|\, u(x)\,
v(x)\, dx&\leq& C_\ell\int_{\zR^n}{\cal
M}_{{\Phi}_1,L}\vec{f}(x)\, M_{L(\log L)^{\ell}}u(x)\, v(x)\, dx\\
&& + \ell C_\ell\int_{\zR^n}{\cal M}_{{\Phi}_1,
L\log L}\vec{f}(x)\, Mu(x)\, v(x)\, dx.
\end{eqnarray*}
\end{corollary}

\medskip

\Proof : Since $v\in RH_{\infty}$ then $v\in A_{\infty}$. Thus, if
$|Q|\leq C|E|$ then $v(Q)\leq C v(E)$. From Lemma $\ref{discretizacion}$ with $q=1$
we have that
\begin{eqnarray*}
\int_{\zR^n}|{\cal T}_{b_j^\ell,\phi}(\vec{f})| \ u \ v &\leq&
C_\ell\sum_{k,\eta}{\Phi}_1( l(Q_{k,\eta}^{0})) ||u ||_{L(\log
L)^\ell, 3Q_{k,\eta}^{0}}
\prod_{i=1}^m ||f_i||_{L,3 Q_{k,\eta}^0} \sup_{3Q_{k,\eta}^0} v \ |E_{k,\eta}^0|\\
& & + \ell C_\ell \sum_{k,\eta}{\Phi}_1(
l(Q_{k,\eta}^{j})) ||u ||_{L,3 Q_{k,\eta}^{j}}
\prod_{i=1}^m ||f_i||_{L(\log L)^{\delta_{i,j}},3
Q_{k,\eta}^{j}}\sup_{3Q_{k,\eta}^{j}}v\ |E_{k,\eta}^{j}|\\
&\leq& C_\ell\sum_{k,\eta}{\Phi}_1(
l(Q_{k,\eta}^{0})) ||u ||_{L(\log L)^\ell, 3Q_{k,\eta}^{0}}
\prod_{i=1}^m ||f_i||_{L,3 Q_{k,\eta}^0} v(E_{k,\ell}^0)\\
& & + \ell C_\ell\sum_{k,\eta}{\Phi}_1(
l(Q_{k,\eta}^{j})) ||u ||_{L,3 Q_{k,\eta}^{j}}
\prod_{i=1}^m ||f_i||_{L(\log L)^{\delta_{i,j}},3
Q_{k,\eta}^{j}}v(E_{k,\eta}^{j})\\
&\leq& C_\ell\sum_{k,\eta}\int_{E_{k,\eta}^0}{\cal
M}_{{\Phi}_1,L}(\vec{f})\ M_{L(\log L)^\ell}(u)\ v\\
& & + \ell C_\ell\sum_{k,\eta}\int_{E_{k,\eta}^{j}}{\cal
M}_{{\Phi}_1, L(\log L)}(\vec{f})\ M(u)\ v.\\
\end{eqnarray*}

\bigskip


\begin{lemma}\label{RH}
Let $g$ be a function such that $M(g)$ is finite a.e. and
$\alpha>0$. Then $M(g)^{-\alpha}\in RH_{\infty}$.
\end{lemma}


\vspace{0.8cm}

\noindent {\it Proof of theorem $\ref{for-t-d}$} :

\medskip

\noindent i) We follow similar arguments to those in \cite{CPSS},
(see also \cite{P} for the multilinear case). We shall use the
duality in $L^p(\mu)$ spaces for $p<1$: if $f\ge0$
\begin{equation*}
\|f\|_{L^p(\mu)}=\inf\{\int fu^{-1} \,d\mu:\|u^{-1}\|_{L^{p'}(\mu)}=1\}=\int fu_0^{-1}\, d\mu,
\end{equation*}
for some $u_0\ge0$ such that $\|u_0^{-1}\|_{L^{p'}(\mu)}=1$, with
$p'=\frac{p}{p-1}<0$. This follows from the following reverse
H\"{o}lder's inequality, which is a consequence of the
H\"{o}lder's inequality,
\begin{equation}\label{rh}
\int fg \, d\mu\ge \|f\|_{L^p(\mu)} \|g\|_{L^{p'}(\mu)}.
\end{equation}

\noindent Given a weight $u$, we shall use the above results with
an absolutely continuous measure $\mu$ with density $M_{L(\log
L)^\ell}u$. In fact, let $g$ be a nonnegative function such that
$\|g^{-1}\|_{L^{p'}(M_{L(\log L)^\ell}u)}=1$ and
\begin{equation*}
\|{\cal M}_{{\Phi}_1,
L(\log L)^\ell}\vec{f}\|_{L^{p}(M_{L(\log L)^\ell}u)}=\int {\cal
M}_{{\Phi}_1, L(\log L)^\ell}\vec{f}\ \frac{M_{L(\log L)^\ell}u}{g}.
\end{equation*}
Let $\delta >0$. By the Lebesgue differentiation theorem we get
\begin{equation*}
\|{\cal M}_{{\Phi}_1,
L(\log L)^\ell}\vec{f}\|_{L^{p}(M_{L(\log L)^\ell}u)}\ge \int {\cal
M}_{{\Phi}_1, L(\log L)^\ell}\vec{f}\
\frac{M_{L(\log L)^\ell}u}{M_{\delta}(g)},
\end{equation*}
where $M_{\delta}(g)=M(g^{\delta})^{1/\delta}$. Then applying
Lemma $\ref{RH}$ and Corollary $\ref{TM}$ to the weight
$M_{\delta}(g)^{-1}$ and the reverse H\"{o}lder's inequality
$(\ref{rh})$, we obtain that

\begin{eqnarray*}
(C_{\ell}+\ell C_\ell)\|{\cal
M}_{{\Phi}_1, L(\log L)^\ell}\vec{f}\|_{L^{p}(M_{L(\log L)^\ell}u)}
&\ge&C_\ell\int {\cal
M}_{{\Phi}_1,L}\vec{f}\
\frac{M_{L(\log L)^\ell}u}{M_{\delta}(g)}\\
&&\quad+\ell C_\ell \int {\cal M}_{{\Phi}_1,
L(\log L)^\ell}\vec{f}\
\frac{Mu}{M_{\delta}(g)}\\&\ge&\int_{\zR^n}|{\cal
T}_{{\vec{b}}^\ell,\phi}\vec{f}(x)|\, \frac{u(x)}{M_{\delta}(g)}\, dx\\
&\ge&\|{\cal
T}_{{\vec{b}}^\ell,\phi}\vec{f}\|_{L^p(u)}\|M_{\delta}(g)^{-1}\|_{L^{p'}(u)}.
\end{eqnarray*}
To finish the proof we shall show that
\begin{equation*}
\|M_{\delta}(g)^{-1}\|_{L^{p'}(u)}\ge\|g^{-1}\|_{L^{p'}(Mu)}=1.
\end{equation*}
Since $p'<0$, this is equivalent to prove that
\begin{equation*}
\int_{\zR^n}M_{\delta}(g)^{-p'}(x) u(x)\, dx\leq C\int
g^{-p'}(x)Mu(x)\, dx.
\end{equation*}
By choosing $\delta$ such that $0<\delta<\frac{p}{1-p}$, we have
that $-p'/\delta>1$ and the above inequality follows from the
classical weighted norm inequality of Fefferman-Stein (see
\cite{FS}). \\

\noindent ii) It is enough to
prove that
\begin{eqnarray*}
\int_{\zR^n}|{\cal
T}_{{\vec{b}}^\ell,\phi}(\vec{f})(x)|u(x)^{1/p}g(x)\, dx&\leq&
C_\ell \left(\int_{\zR^n}[{\cal M}_{{\Phi}_1, L(\log
L)^\ell}\vec{f}(x)]^p M_{L(\log L)^{[\ell p+p]}}u(x)\,
dx\right)^{1/p}
\|g\|_{p'}\\
\end{eqnarray*}
for every function $g\in L^{p'}(\zR^n)$. From Lemma $\ref{discretizacion}$ we get that
\begin{eqnarray*}\label{discreta1}
\int_{\zR^n}|{\cal T}_{b_j^\ell,\phi}(\vec{f})| \ u^{1/p} \ g
&\leq& C_\ell\sum_{k,\eta}{\Phi}_1( l(Q_{k,\eta}^{0})) ||u^{1/p}
g||_{L(\log L)^\ell, 3Q_{k,\eta}^{0}}
 \prod_{i=1}^m ||f_i||_{L^1,3 Q_{k,\eta}^0} |E_{k,\eta}^0|\nonumber\\
 & & + \ell C_\ell \sum_{k,\eta}{\Phi}_1(l(Q_{k,\eta}^{j}))
 ||u^{1/p} g||_{L,3 Q_{k,\eta}^{j}}
 \prod_{i=1}^m ||f_i||_{L(\log L)^{\delta_{i,j}},3 Q_{k,\eta}^{j}}
 |E_{k,\eta}^{j}|.
\end{eqnarray*}
Let $\epsilon_0$ and $\epsilon_1$ be two positive numbers to be
chosen later and consider the Young functions
$\mathcal{A}_{\ell}(t)=t^p(1+\log^+t)^{\ell p+(p-1)(1+\epsilon_{\ell})}$ and
$\mathcal{C}_{\ell}(t)=t^{p'}(1+\log^+t)^{-(1+\epsilon_{\ell})}$. Then, it
is easy to check that $\mathcal{A}_{\ell}^{-1}(t)\mathcal{C}_{\ell}^{-1}(t)\leq
\mathcal{B}_{\ell}^{-1}(t)$, where $\mathcal{B}_\ell(t)=t(1+\log^+ t)^\ell$, and $\mathcal{C}_{\ell}\in B_{p'}$. Thus, by the generalized
H\"{o}lder's inequality
\begin{equation*}
M_{L(\log L)^{\ell}}(u^{1/p}g)\leq M_{{\mathcal{A}_\ell}}(u^{1/p})M_{{\mathcal{C}_\ell}}g=
\left(M_{{\tilde{\mathcal{A}}_{\ell}}}u\right)^{1/p} \, M_{{\mathcal{C}_{\ell}}}g
\end{equation*}
where $\tilde{\mathcal{A}}_{\ell}(t)=t(1+\log^+ t)^{\ell
p+(p-1)(1+\epsilon_{\ell})}$. Then
\begin{eqnarray*}
\int_{\zR^n}|{\cal T}_{b_j^\ell,\phi}(\vec{f})| u^{1/p}  g &\leq &
C_{\ell} \int_{\zR^n}{\cal M}_{{\Phi}_1,L}\vec{f}(x)
M_{{\mathcal{A}}_{\ell}}(u^{1/p})(x)
M_{\mathcal{C}_{\ell}}(g)(x)\, dx\\
 & & + \ell C_\ell \int_{\zR^n}{\cal
M}_{{\Phi}_1,L\log L}\vec{f}(x)
M_{{\mathcal{A}}_{0}}(u^{1/p})(x)
M_{\mathcal{C}_{0}}(g)(x)\, dx\\
&\leq&C_\ell \int_{\zR^n}{\cal
M}_{{\Phi}_1,L(\log L)^{\ell}}\vec{f}
M_{{\mathcal{A}}_{\ell}}(u^{1/p})
[M_{\mathcal{C}_{0}}g+M_{\mathcal{C}_{1}}g]\\
&\leq&C_\ell\left(\int_{\zR^n}[{\cal
M}_{{\Phi}_1,L(\log L)^{\ell}}\vec{f}
M_{{\mathcal{A}}_{\ell}}(u^{1/p})]^{p}\right)^{1/p}
\|g\|_{p'}\\
&\leq&C_\ell\left(\int_{\zR^n}[{\cal
M}_{{\Phi}_1,L(\log L)^{\ell}}\vec{f}]^p
M_{\tilde{\mathcal{A}}_{\ell}}(u)\right)^{1/p} \|g\|_{p'}
\end{eqnarray*}

Thus, by taking $\epsilon_{\ell}=\frac{[\ell p+p]-\ell
p-p+1}{p-1}$ we are done. $\square$

\bigskip

\section{Proofs of the weak type results}

\noindent {\it Proof of theorem $\ref{debilmaximal}$} : Let $\Omega_{\lambda}=\{x\in \zR^n: {\cal M}_{\varphi,\mathcal
B}\vec{f}(x)>\lambda^m\}$. By homogeneity we may assume that
$\lambda=1$. Let $K$ be a compact set in $\Omega_{1}$. Since $K$
is a compact set and using Vitali's covering lemma we obtain a
finite family of disjoint cubes $\{Q_j\}$ for which
\begin{equation*}\label{cubos}
\varphi(|Q_j|)\prod_{i=1}^m \|f_i\|_{\mathcal B,Q_j}>1,
\end{equation*}
or, equivalently
\begin{equation*}
\prod_{i=1}^m \|\varphi(|Q_j|)^{1/m}f_i\|_{\mathcal B,Q_j}>1
\end{equation*}
and $K\subset U_j3Q_j$. The proof follows now similar arguments as in the proof of Theorem
4.1 of \cite{PPTT}. We include it for the sake of completeness.

Let $g_i=\varphi(|Q_j|)^{1/m}f_i$ and $C_h^m$ the family
of all subset $\sigma=\{\sigma(1),...,\sigma(h)\}$ of different elements from the index
$\{1,...,m\}$ with $1\leq h\leq m$ . Given
$\sigma\in C_h^m$ and a cube $Q_i$, we say that $i\in B_{\sigma}$
if $\|g_{\sigma(k)}\|_{\mathcal B,Q_i}>1$ for $k=1,...,h$ and
$\|g_{\sigma(k)}\|_{\mathcal B,Q_i}\leq1$ for $k=h+1,...,m$. For  $\sigma\in C_h^m$ and $i\in B_{\sigma}$ denote
$$\Pi_k=\prod_{j=1}^k\|g_{\sigma(j)}\|_{\mathcal B, Q_i}$$
and $\Pi_0=1$. Then  $\Pi_k>1$ for every $1\leq k\leq m$ and thus
\begin{equation*}
1<\Pi_{k}=\|g_{\sigma(k)}\|_{\mathcal B, Q_i}\,
\Pi_{k-1}=\|g_{\sigma(k)}\Pi_{k-1}\|_{\mathcal B, Q_i}
\end{equation*}
or, equivalently
\begin{equation}\label{12}
\frac{1}{|Q_i|}\int_{Q_i} \mathcal B\left(g_{\sigma(k)}\,
\Pi_{k-1}\right)>1.
\end{equation}
In particular,
\begin{equation}\label{13}
1<\frac{1}{|Q_i|}\int_{Q_i} \mathcal B\left(g_{\sigma(m)}\,
\Pi_{m-1}\right)\leq\frac{1}{|Q_i|}\int_{Q_i} \mathcal
B\left(g_{\sigma(m)}\right) \mathcal B\left(\Pi_{m-1}\right).
\end{equation}
Now, by taking into account the equivalence
\begin{equation*}
\|g\|_{{\mathcal B}, Q}\simeq
\inf_{\mu>0}\left\{\mu+\frac{\mu}{|Q|}\int_Q \mathcal B(|g|/\mu)\right\},
\end{equation*}
if $1\leq j\leq m-h-1$, by $(\ref{12})$ we get
\begin{eqnarray*}
\mathcal B^j(\Pi_{m-j})&=&\mathcal
B^j(\|g_{\sigma(m-j)}\Pi_{m-j-1}\|_{{\mathcal B},
Q_i})\\
&\leq&C\mathcal B^j\left(1+\frac{1}{|Q_i|}\int_{Q_i}
\mathcal B\left(g_{\sigma(m-j)}\, \Pi_{m-j-1}\right)\right)\\
&\leq&\frac{C}{|Q_i|}\int_{Q_i} \mathcal
B^{j+1}\left(g_{\sigma(m-j)}\right)\, \mathcal
B^{j+1}\left(\Pi_{m-j-1}\right).
\end{eqnarray*}
From $(\ref{13})$, by iterating the inequality above, we obtain
\begin{eqnarray*}
1&<&C\frac{1}{|Q_i|}\int_{Q_i} \mathcal
B\left(g_{\sigma(m)}\right)\frac{1}{|Q_i|}\int_{Q_i}
\mathcal B^2\left(g_{\sigma(m-1)}\right)\mathcal B^2\left(\Pi_{m-2}\right)\\
&\leq&C\left(\prod_{j=0}^{m-h-1}\frac{1}{|Q_i|}\int_{Q_i}
\mathcal B^{j+1}\left(g_{\sigma(m-j)}\right)\right)\mathcal B^{m-h}\left(\Pi_{h}\right)\\
&\leq&C \left(\prod_{j=0}^{m-h-1}\frac{1}{|Q_i|}\int_{Q_i}\mathcal
B^{j+1}\left(g_{\sigma(m-j)}\right)\right)
\left(\prod_{j=1}^h\mathcal B^{m-h}(\|g_{\sigma(j)}\|_{{\mathcal
B}, Q_i})\right).
\end{eqnarray*}
since $\mathcal B$ is submultiplicative. Thus, since $i\in B_{\sigma}$, we have
$\|g_{\sigma(j)}\|_{{\mathcal B}, Q_i}>1$ for $j=1,...h$, and it
follows that
\begin{equation}\label{cuatro}
1<C\left(\prod_{j=0}^{m-h-1}\frac{1}{|Q_i|}\int_{Q_i}\mathcal
B^{j+1}\left(g_{\sigma(m-j)}\right)\right)
\left(\prod_{j=1}^h\frac{1}{|Q_i|}\int_{Q_i}\mathcal
B^{m-h+1}(g_{\sigma(j)})\right).
\end{equation}
Now, since for $1\leq h\leq m$ and $0\leq j\leq m-h-1$ we have
that $\mathcal B^{j+1}(t)\leq \mathcal B^{m-h}(t)\leq \mathcal
B^m(t)$ and $\mathcal B^{m-h+1}(t)\leq \mathcal B^m(t)$. Then
\begin{equation*}
1<C\prod_{j=1}^m\frac{1}{|Q_i|}\int_{Q_i}\mathcal B^{m}(g_{j})
\end{equation*}
or equivalently
\begin{equation*}
|Q_i|<C\prod_{j=1}^m\left(\int_{Q_i}\mathcal
B^{m}(g_j)\right)^{1/m}.
\end{equation*}
Thus, we obtain that
\begin{equation*}\label{nueva}
|Q_{j}|\leq C\prod_{i=1}^{m}\left(\int_{Q_j}\mathcal
B^m\left(\varphi(|Q_j|)^{1/m}|f_i|\right)\right)^{1/m}\leq
C\prod_{i=1}^{m}\left(\mathcal
B^m\left(\varphi(|Q_j|)^{1/m}\right)\int_{Q_j}\mathcal
B^m\left(|f_i|\right)\right)^{1/m}
\end{equation*}
since $B$ is submultiplicative.

By applying H\"{o}lder's inequality we have that
$\frac{u(Q)}{|Q|}\leq
\prod_{i=1}^m\left(\frac{u_i(Q)}{|Q|}\right)^{1/m}$. Then from the above inequalities we obtain that
\begin{eqnarray*}
u(K)&\leq&C \sum_j
\frac{u(3Q_j)}{|3Q_j|}|Q_j|\\
&\leq&C
\sum_j\prod_{i=1}^m\left(\frac{1}{|3Q_j|}\int_{3Q_j}u_i\right)^{1/m}
\left(\mathcal B^m\left(\varphi(|Q_j|)^{1/m}\right)\int_{Q_j}\mathcal B^m\left(|f_i|\right)\right)^{1/m}\\
&\leq&C \prod_{i=1}^m\left(\sum_j\mathcal
B^m\left(\varphi(|Q_j|)^{1/m}\right)\frac{1}{|3Q_j|}\int_{3Q_j}u_i
\int_{Q_j}\mathcal B^m\left(|f_i|\right)\right)^{1/m}\\
&\leq&C \prod_{i=1}^m\left(\int_{\zR^n}\mathcal B^m\left(|f_i|\right)M_{\psi,L}u_i\right)^{1/m}\\
\end{eqnarray*}
where $\psi=\mathcal B^m\circ \varphi^{1/m}$ and the proof
concludes. $\square$

\bigskip

In order to prove Theorem \ref{control} we need the following lemma.

\begin{lemma} Let $q>1$, $\epsilon>0$, $\mathcal B_{\ell}(t)=t(1+\log^+t)^{\ell}$,
$\tilde{\mathcal{A}}_{\ell}(t)=t(1+\log^+t)^{\ell q+q-1+\epsilon}$ and
let $w$, $u$ be weights. Then there exists a positive constant $C$
such that the inequality
\begin{equation}\label{dos}
\int_{\zR^n}M_{\mathcal
B_{\ell}}(f)^{q'}M_{\tilde{\mathcal{A}}_{\ell}}(u)^{1-q'}w\leq
\int_{\zR^n}|f|^{q'}u^{1-q'}Mw
\end{equation}
holds.
\end{lemma}

\noindent \Proof : Let ${\mathcal{A}}_{\ell}(t)=t^q(1+\log^+t)^{\ell
q+q-1+\epsilon}$ and
$\mathcal{C}(t)=t^{q'}(1+\log^+t)^{-(1+\frac{\epsilon}{q-1})}$. Then, it is
easy to check that $\mathcal{A}_{\ell}^{-1}(t)\mathcal{C}^{-1}(t)\leq \mathcal
B_{\ell}^{-1}$(t). Thus, by H\"{o}lder's inequality
$$M_{\mathcal{B}_{\ell}}(gu^{1/q})\leq
M_{\mathcal{C}}(g)M_{\mathcal{A}_{\ell}}(u^{1/q})=M_{\mathcal{C}}(g)
M_{\tilde{\mathcal{A}}_{\ell}}(u)^{1/q}.$$
Then, since $\mathcal{C}\in B_{q'}$
\begin{eqnarray*}
\int_{\zR^n}M_{\mathcal
B_{\ell}}(gu^{1/q})^{q'}M_{\tilde{\mathcal{A}}_{\ell}}(u)^{1-q'}w&\leq&
C\int_{\zR^n}M_{\mathcal{C}}(g)^{q'}w
\leq C\int_{\zR^n}|g|^{q'}Mw,
\end{eqnarray*}
where the last inequality was proved in \cite{Pe2}.
Finally, by taking $g=f u^{-1/q}$ we obtain ($\ref{dos}$). $\square$

\bigskip

\noindent {\it Proof of theorem $\ref{control}$}: Let $p>1$ to be
chosen later. Thus, since $L^{p,\infty}$ and $L^{p',1}$ are
associated spaces, we have that
\begin{equation*}
\|{\cal
T}_{\vec{b}^\ell,\phi}\vec{f}\|_{L^{1/m,\infty}(u)}^{1/(pm)}=\||{\cal
T}_{\vec{b}^\ell,\phi}\vec{f}|^{1/(pm)}\|_{L^{p,\infty}(u)}=\sup_{\|g\|_{L^{p',1}(u)}\leq
1}\int_{\zR^n}|{\cal T}_{\vec{b}^\ell,\phi}\vec{f}|^{1/(pm)} g
\, u.
\end{equation*}
By Theorem $\ref{for-t-d}$ i) we obtain that
\begin{eqnarray*}
\int_{\zR^n}|{\cal T}_{\vec{b}^\ell,\phi}\vec{f}|^{1/(pm)} g \,
u&\leq& C_{\ell}\int_{\zR^n}({\cal M}_{{\Phi}_1,
B_{\ell}}\vec{f})^{1/(pm)} M_{B_{\ell}}(gu)\\
&=& C_{\ell}\int_{\zR^n}({\cal M}_{{\Phi}_1,
B_{\ell}}\vec{f})^{1/(pm)}
\frac{M_{B_{\ell}}(gu)}{M_{\tilde{A}_{\ell}}(u)}M_{\tilde{A}_{\ell}}(u),
\end{eqnarray*}
where $\tilde{A}_{\ell}$ is a Young function to be chosen later. Then, by applying
H\"{o}lder's inequality in Lorentz spaces we obtain that
\begin{eqnarray*}
\int_{\zR^n}|{\cal T}_{\vec{b}^\ell,\phi}\vec{f}|^{1/(pm)} g \,
u&\leq&\|({\cal M}_{{\Phi}_1,
B_{\ell}}\vec{f})^{1/(pm)}\|_{L^{p,\infty}(M_{\tilde{A}_{\ell}}(u))}\left\|\frac{M_{B_{\ell}}(gu)}{M_{\tilde{A}_{\ell}}(u)}
\right\|_{L^{p',1}(M_{\tilde{A}_{\ell}}(u))}.
\end{eqnarray*}
Now the argument follows in similar way to that in \cite{CPSS}
(see also \cite{P} for the multilinear fractional integral).
If we define
$
S(f)=\frac{M_{B_{\ell}}(fu)}{M_{\tilde{A}_{\ell}}(u)},
$, the result will be done if we prove that
$S:L^{p',1}(u)\rightarrow L^{p',1}(M_{\tilde{A}_{\ell}}(u))$.

We shall take $\tilde{A}_{\ell}$ such that
$M_{B_{\ell}}(u)\leq M_{\tilde{A}_{\ell}}(u)$, then $S:L^{\infty}(u)\rightarrow
L^{\infty}(M_{\tilde{A}_{\ell}}(u))$. Thus, by the
Marcinkiewicz's interpolation theorem in Lorentz spaces due to
Hunt (see \cite{BS}), it is enough to prove that
$S:L^{(p+\epsilon)'}(u)\rightarrow
L^{(p+\epsilon)'}(M_{\tilde{A}_{\ell}}(u))$, which is equivalent
to prove the following inequality
\begin{equation}\label{dual}
\int_{\zR^n}(M_{B_{\ell}}(f))^{(p+\epsilon)'}(M_{\tilde{A}_{\ell}}(u))^{1-(p+\epsilon)'}\leq
\int_{\zR^n}|f|^{(p+\epsilon)'}u^{1-(p+\epsilon)'}.
\end{equation}
By choosing $\tilde{A}_{\ell}=t(1+\log^+t)^{\ell
p+p-1+(\ell+2)\epsilon}$, $\epsilon>0$, the above result
holds from Lemma $\ref{dual}$ with $q=p+\epsilon$ and $w=1$.
Then, we obtain that
\begin{equation*}
\|{\cal T}_{\vec{b}^\ell,\phi}\vec{f}\|_{L^{1/m,\infty}(u)}\leq
C \|({\cal M}_{{\Phi}_1,
B_{\ell}}\vec{f})\|_{L^{1/m,\infty}(M_{\tilde{A}_{\ell}}(u))}.
\end{equation*}
Thus, taking
$p=1+\frac{\delta_{\ell}-\epsilon(\ell+2)}{\ell+1}$ with
$0<\epsilon(\ell+2)<\delta_{\ell}$ we get the theorem, since $\tilde{A}_{\ell}(t)=t(1+\log^+t)^{\ell+\delta_{\ell}}$.


\begin{thebibliography}{AAAAAA}
\bibitem[BS]{BS} Bennett, C. and Sharpley, R.: \emph{Interpolation of
operators}, Pure and Applied Mathematics, vol. 129, Academic Press
Inc., Boston, MA, 1988.


\bibitem[CN]{CN} Cruz Uribe, D. and Neugebauer, C. J.: \emph{The structure of the reverse H\"{o}lder's classes},
Trans. Amer. Math. Soc. 347 (1995), 2941--2960.

\bibitem[CPSS]{CPSS} Carro, M. J., P\'erez, C., Soria, F. and
Soria, J.: \emph{Maximal functions and the control of weighted
inequalities for the fractional integral operator}, Indiana Univ.
Math. J.  54 (2005),  627-644.

\bibitem[CX]{CX} Chen, X. and Xue, Q.:\emph{Weighted estimates for a class of multilinearfractional operators},
 J. Math. Anal. Appl. 362 (2010), 355-373.


\bibitem[CUP]{CUP} Cruz Uribe, D. and P\'erez, C.: \emph{On the two-weight problem for singular integral operators}, Ann. Scuola Norm. Sup. Pisa CI, Sci, 1 (2002) 821--849.

\bibitem[CUMP]{CUMP} Cruz Uribe, D., Martell, J. M. and P\'erez, C.: \emph{Extrapolation from $A_{\infty}$
weights and applications}, J. Funct. Anal. 213 No. 2 (2004),
412--439.

\bibitem[FS]{FS} Fefferman, C. and Stein, E., M.: \emph{Some maximal
inequalities}, Amer. J. Math. 93 (1971), 107--115.

\bibitem[G]{G} Grafakos, L.: \emph{Classical and Modern Fourier
Analysis}, Pearson Education Inc. (2004).

\bibitem[L1]{L1} Li, W.: \emph{Weighted inequalities for commutators of potential type operators}, J. Korean Math. Soc. 44  No. 6 (2007), 1233--1241.

\bibitem[L2]{L2} Li, W.: \emph{Two-weight norm inequalities for commutators of potential type integral operators}, J. Math. Anal. Appl. 322 (2006), 1215--1223.

\bibitem[LOPTT]{LOPTT} Lerner, A., Ombrosi, S., P\'erez, C., Torres, R., and Trujillo-Gonz\'alez, R.:
\emph{ New maximal functions and multiple weights for the
multilinear Calder\'on-Zygmund theory} Adv. Math. {\bf 220}
(4)(2009), 1222-1264

\bibitem[LQY]{LQY} Li, W. M., Qi, J. Y. and Yan, X. F.: \emph{Weighted norm   inequalities for potential type operators}, J. Math. Research Exposition 29 (2009), 895--900.

\bibitem[LYY]{LYY} Li, W., Yan, X. and Yu, X.: \emph{Two-weight inequalities for commutators of potential operators on spaces of homogeneous type}, Potential Anal. 31 (2009), 117--131.

    \bibitem[M]{M} Moen, K.:  \emph{Weighted inequalities for multilinear fractional integral operators},
Collect. Math. 60 (2) (2009), 213-238


\bibitem[MMN]{MMN} Maldonado, D., Moen, K. and Naibo, V.:
 \emph{Weighted multilinear Poincar`\'e inequalities for vector fields of H\"ormander
 type}, to appear in Indiana Univ. Math. J.

\bibitem[MW]{MW} Muckenhoupt, B. and Wheeden R.:, \emph{Weighted norm inequalities for fractional
integrals}, Trans. Amer. Math. Soc. 192 (1974) 261–274.

\bibitem[O]{O} O'Neil, R.: \emph{Fractional integration in Orlicz spaces. I.}, Trans. Amer. Math. Soc. 115 (1965), 300--328.

\bibitem[Pe]{Pe} P\'erez, C.: \emph{Two weighted inequalities for potential and fractional type maximal operators},
Indiana Univ. Math. 43 (2) (1994) 663--683.
\bibitem[Pe1]{Pe1} P\'erez, C.: \emph{Sharp $L^p$-weighed Sobolev
inequalities},Ann.de l'inst. Fourier 45 (3) (1995) 809--824.
\bibitem[Pe2]{Pe2}P\'erez, C.: \emph{On sufficient conditions for the
boundedness of the Hardy-Littlewood maximal operator between
weighted $L^p$-spaces with different weights }, Proc. London Math.
Soc. (3) 71 (1995) 135--157.

\bibitem[Pe3]{Pe3} P\'erez, C.:\emph{Endpoint estimates for commutators of singular integral operators},
 J. Funct. Anal. {\bf 128} (1995), 163--185.

\bibitem[PW]{PW} P\'erez, C. and Wheeden, R.: \emph{Uncertainty principle estimates for vector
fields}, J. Funct. Anal. 181, (2001), 146--188.


\bibitem[PPTT]{PPTT} P\'erez, C., Pradolini, G., Torres, R. H. and
Trujillo-Gonz\'alez, R. \emph{Endpoint estimates for iterated
commutators of multilinear singular integrals},
preprint.

\bibitem[P]{P} Pradolini, G.: \emph{Weighted inequalities
and pointwise estimates for the multilinear fractional integral
and maximal operators}, J. Math. Anal. and Appl. 367, (2), (2010),
640-656

\bibitem[SW]{SW} Sawyer, E. and Wheeden, R.:\emph{Weighted
inequalities for fractional integrals on euclidean and homogeneous
spaces}, Amer. J. Math. 114 (1992), 813-874.




\end{thebibliography}
\end{document}